\documentclass[10pt,twoside]{article}
\usepackage[latin1]{inputenc}
\usepackage{amsmath}
\usepackage{ epsfig}
\usepackage{graphicx}

\input epsf
\def\limiten{\renewcommand{\arraystretch}{0.5}
\begin{array}[t]{c}\stackrel{}{\longrightarrow} \\
{\scriptstyle n\rightarrow
\infty}\end{array}\renewcommand{\arraystretch}{1}}

\def\limitepsn{\renewcommand{\arraystretch}{0.5}
\begin{array}[t]{c}\stackrel{a.s.}{\longrightarrow} \\
{\scriptstyle n \rightarrow
\infty}\end{array}\renewcommand{\arraystretch}{1}}

\def\limiteloin{\renewcommand{\arraystretch}{0.5}
\begin{array}[t]{c}\stackrel{{\cal D}}{\longrightarrow} \\
{\scriptstyle n\rightarrow
\infty}\end{array}\renewcommand{\arraystretch}{1}}

\newcommand{\be}{\begin{equation}}
\newcommand{\ee}{\end{equation}}
\newcommand{\bd}{\begin{displaymath}}
\newcommand{\ed}{\end{displaymath}}

\newcommand{\ba}{\begin{eqnarray}}
\newcommand{\ea}{\end{eqnarray}}
\newcommand{\ban}{\begin{eqnarray*}}
\newcommand{\ean}{\end{eqnarray*}}

\newcommand{\R} {I\!\!R}
\newcommand{\E} {E\,}
\newcommand{\N} {I\!\! N}
\newcommand{\Z} {{\bf Z}}

\renewcommand{\arraystretch}{.8}

\renewcommand{\Box}{\hfill\rule{0.25cm}{0.25cm}} 

\newtheorem{lem}{Lemma}[section]
\newtheorem{Theo}{Theorem}[section]
\newtheorem{cor}{Corollary}[section]
\newtheorem{rem}{Remark}[section]
\newenvironment{dem}{\ \\ {\bf Proof. }}
{\Box\par\medskip\noindent}

\def\1{{\bf 1}}

\begin{document}
\title{\bf Testing for parameter constancy in general causal time series models}
 \maketitle \vspace{-1.0cm}
\begin{center}
 by William~Charky~KENGNE \footnote{ Supported by AUF (Agence Universitaire de la Francophonie) and Edulink ACP-EU project.}
\end{center}
\begin{center}
 {\it SAMM, Université Paris 1 Panth\'eon-Sorbonne, 90 rue de Tolbiac 75634-Paris Cedex 13, France.\\
  E-mail: William-Charky.Kengne@malix.univ-paris1.fr}
\end{center}

 \pagestyle{myheadings}
 \markboth{Testing for parameter constancy in general causal time series models}{ William KENGNE}

\textbf{Abstract} :  We consider a process $ X= (X_t)_{t\in \Z}$
belonging to a large class of
   causal models including AR($\infty$), ARCH($\infty$), TARCH($\infty$),... models.
   We assume that the model depends on a parameter $\theta_0 \in \R^d$ and consider the problem of testing for change
   in the parameter. Two statistics $\widehat{Q}^{(1)}_n$ and $ \widehat{Q}^{(2)}_n$ are constructed using quasi-likelihood
   estimator (QLME) of the parameter. Under the null hypothesis that there is no change, it is shown that each of these two statistics weakly converges
   to the supremum of the sum of the squares of independent Brownian bridges.
    Under the local alternative that there is one change, we show that the test statistic $ \widehat{Q}_n=\text{max} \big(\widehat{Q}^{(1)}_n , \widehat{Q}^{(2)}_n  \big) $
    diverges to infinity.
   Some simulation results for AR(1), ARCH(1) and GARCH(1,1) models are reported to show the applicability and the performance of our
   procedure with comparisons to some other approaches.~\\ ~\\


\section{Introduction}
   Many statistical data can be represented by models which may change over
 time, for instance hydraulic flow, climate data. Before any inference on these data, it is crucial to test whether  a change
  has not occurred in the model.\\
  \indent Since Page \cite{Page1955} in 1955, real advances have been done about tests for change detection.
  Horvath \cite{Horvart1993} proposed a test for detecting a change in the parameter of
  autoregressive processes  based on weighted supremum and $L_p$-functionals of the residual sums.
   The CUSUM statistic
   which  was introduced by Brown {\it et al.} \cite{Brown1975} in 1975, was modified by Inclan and Tiao \cite{Tiao1994}
   for testing change in variance of independent random variables.
   Their test has asymptotically correct size but the asymptotic power is unknown.
   Numerous works devoted to the CUSUM-type procedure, for instance Kim  {\it et al.} \cite{Kim2000}  for testing change in parameters of GARCH(1,1),
    Kokoszka and Leipus \cite{Leipus1999} in the specific case of ARCH($\infty$) or Aue {\it et al.} \cite{Horman2009} for testing breaks in
   covariance.  Kulperger and Yu \cite{Kulperger2005} studied the high moment
  partial sum process based on residuals  and applied it to residual CUSUM test in GARCH model.
  Horváth {\it et al.} \cite{Horvart2008}  suggested to compute the ratio of the CUSUM functionals instead of the
 differences for testing change in the mean of a time series. Berkes
 {\it et al.} \cite{Berkes2004} used a test based on approximate likelihood scores
 for testing parameter constancy in GARCH(p,q) models.  These procedures are mostly developed in a parametric framework and their asymptotic powers are unknown.
 The present work is a new contribution to the challenging problem of  test for change detection.\\

 In this paper, we consider a general class ${\cal M}_T(M,f)$ of causal (non-anticipative) time series. Let $M , f : \R^{\N} \rightarrow \R$ be
 measurable functions,  $(\xi_t)_{t\in \Z}$ be a sequence of centered independent
 and identically distributed (iid)  random variables called the innovations and satisfying $\textrm{var}(\xi_0)= \sigma^2$  and $\Theta $ a compact subset of $\R^{d}$.
 Let $T \subset \Z$, and for any $ \theta \in \Theta$, define ~\\

{\bf Class ${\cal M}_T(M_{\theta},f_{\theta})$:} {\it The process
$X=(X_t)_{t\in\Z}$ belongs to  ${\cal M}_T(M_{\theta},f_{\theta})$
if it satisfies the relation:}
 \begin{equation}\label{model}
 X_{t+1}=M_{\theta} \big ((X_{t-i})_{i\in \N}\big ) \xi_t+f_{\theta}\big
 ((X_{t-i})_{i\in \N}\big )\quad\mbox{for all $t\in T$}.
 \end{equation}

  \noindent The existence and properties of this general class of affine processes  were studied in Bardet and Wintenberger \cite{Bardet2009}.
  Numerous classical time series are included in ${\cal M}_\Z(M,f)$:
  for instance AR$(\infty)$, ARCH$(\infty)$, TARCH$(\infty)$, ARMA-GARCH or bilinear processes.  \\

   Now, assume that a trajectory $(X_1,\cdots,X_n)$ of $X=(X_t)_{t\in\Z}$ is
   observed and consider the following hypothesis:\\

      \noindent  $H_0$ :   there exists $\theta_0 \in \Theta $ such that $ (X_1, \cdots, X_n)$ belongs to the
       class $\mathcal{M}_{ \{1,\cdots,n \} }(M_{\theta_0},f_{\theta_0})$ ; \\
       $H_1$ : there exist $ K \geq 2 $, $ \theta_1,\cdots,\theta_K \in \Theta  $  such that
       $ (X_1, \cdots, X_n)$ belongs to $ \underset{ j=1} { \overset{ K} { \bigcap } } \mathcal{M}_{T^n_j}(M_{\theta_j},f_{\theta_j})  $
       where $T^n_j=\{t_{j-1}+1,t_{j-1}+2,\cdots,t_j \}$ with $0=t_0<t_1<\cdots<t_{K-1}<t_K=n$ .\\

  Thus, it is easy to see that under $H_1$  the property of stationary  is lost after the first change.
  This is not the case in many existing works (for instance Kouamo {\it et al.} \cite{Kouamo2010} ) where the stationarity or the K-th
  order stationarity after the change is an essential assumption.\\

  In this paper we study a new test for change detection (see Bardet {\it et al.}  \cite{Bardet2010} for the procedure of the estimation of the instants of change).
  We consider a semi-parametric test statistic based on the QLME   which is a modification of  the
  statistic introduced by Lee {\it et al.} \cite{Lee2003}.
  For $k,k'\in \{2,\cdots,n-1 \}$ (with $k\leq k'$) let $\widehat{\theta}_n(X_k,\cdots,X_{k'})$ be the QLME of
  the parameter computed on $\{k,\cdots,k' \}$. The basic idea of our  procedure is that : under $H_0$,
  $\widehat{\theta}_n(X_1,\cdots,X_{k})$ and $\widehat{\theta}_n(X_{k+1},\cdots,X_{n})$ are close to
  $\widehat{\theta}_n(X_1,\cdots,X_{n})$ and the distances  $ \| \widehat{\theta}_n(X_1,\cdots,X_{k}) - \widehat{\theta}_n(X_{1},\cdots,X_{n}) \|$
  and $ \| \widehat{\theta}_n(X_{k+1},\cdots,X_{n}) - \widehat{\theta}_n(X_{1},\cdots,X_{n})\|$ are not too large.
  Therefore, we show that the test statistic is finite under the null hypothesis and diverges to infinity under the alternative of one
  change. Simulation results compared to some other procedures show that our procedure is more powerful.
  In Section 2 we present assumptions, some examples and construct the test
  statistic. In Section 3 we give some asymptotic results. The
  empirical studies of AR(1), ARCH(1) and GARCH(1,1) are detailed in
  Section 4 and the proofs of the main results are presented in Section 5.

\section{Assumptions and test statistics}

\subsection{Assumptions on the class of models ${\cal M}_\Z (f_\theta,M_\theta)$}

 Let $\theta\in\R^d$ and $M_\theta$ and $f_\theta$ be numerical functions such that for all $(x_i)_{i\in \N}\in \R^{\N}$,
 $M_\theta\big((x_i)_{i\in \N}\big)\neq 0$ and $f_\theta\big((x_i)_{i\in \N}\big) \in \R$. We will use the following
 norms:
 \begin{enumerate}
 \item  $\|\cdot\|$ applied to a vector denotes the Euclidean norm of
    the vector;
    \item for any compact set $ \Theta\subseteq\R^d$ and for any
    $g:\Theta \longrightarrow\R^{d'}$, $ \|g\|_\Theta=\sup_{\theta\in \Theta}(\|g(\theta)\|)$.
 \end{enumerate}

\noindent Let  $\Psi_\theta=f_\theta, \, M_\theta$ and $i=0,\, 1, \,
2$, then for any compact set $ \Theta\subseteq\R^d$, define\\
~\\
{\bf Assumption A$_i(\Psi_\theta,\Theta)$}: {\em Assume that
$\|{\partial^i\Psi_\theta(0)}/{\partial\theta^i}\|_\Theta<\infty$
 and there exists a sequence of non-negative real number $(\alpha^{(k)}_i(\Psi_\theta,\Theta))_{i\geq 1}$ such that $ \sum\limits_{k=1}^{\infty} \alpha^{(i)}_k(\Psi_\theta,\Theta) <\infty$  satisfying
\begin{equation*}
\Big\|\dfrac{\partial^i\Psi_\theta(x)}{\partial\theta^i}-\dfrac{\partial^i\Psi_\theta(y)}{\partial\theta^i}\Big\|_\Theta
\leq \sum\limits_{k=1}^{\infty}\alpha^{(i)}_k(\Psi_\theta,\Theta)|x_k-y_k|\quad \mbox{for all}~x, y \in \R^{\N}.\\
\end{equation*}}
 \noindent In the sequel we refer to the particular case called "ARCH-type process"
 if $f_\theta =0$ and if the following assumption holds with $ h_\theta := M_\theta^2$ :\\
~\\
{\bf Assumption A$_i(h_\theta,\Theta)$}: {\em Assume that
$\|{\partial^i h_\theta(0)}/{\partial\theta^i}\|_\Theta<\infty$
 and there exists a sequence of non-negative real number $(\alpha^{(k)}_i(h_\theta,\Theta))_{i\geq 1}$ such as $ \sum\limits_{k=1}^{\infty} \alpha^{(i)}_k(h_\theta,\Theta) <\infty$  satisfying
\begin{equation*}
\Big\|\dfrac{\partial^i
h_\theta(x)}{\partial\theta^i}-\dfrac{\partial^i
h_\theta(y)}{\partial\theta^i}\Big\|_\Theta
\leq \sum\limits_{k=1}^{\infty}\alpha^{(i)}_k(h_\theta,\Theta)|x^2_k-y^2_k|\quad \mbox{for all}~x, y \in \R^{\N}.\\
\end{equation*}}
Then define the set:
   \begin{multline*}
  \Theta(r):=\{\theta \in \Theta,\, A_0(f_\theta,\{\theta\})\, \,\mbox{and} \,A_0(M_\theta,\{\theta\}) \, \textrm{hold with}\,
 \sum\limits_{k \geq 1}^{ } \alpha^{(0)}_k(f_\theta,\theta) + (\E| \xi_0 |^r)^{1/r}\sum\limits_{k \geq 1}^{ } \alpha^{(0)}_k(M_\theta,\theta) <1 \}\\
  \cup \{\theta \in \Theta,\, f_\theta = 0  \textrm{ and}\,  \,A_0(h_\theta,\{\theta\}) \, \textrm{hold with}\,
  (\E| \xi_0 |^r)^{2/r}\sum\limits_{k \geq 1}^{ } \alpha^{(0)}_k(h_\theta,\theta) <1 \}.
 \end{multline*}
 The Lipschitz-type hypothesis $A_i(\Psi_\theta,\Theta)$ are classical when studying the existence of solutions of the general model.
 If $\theta\in\Theta(r)$ the existence of a unique causal, stationary and ergodic  solution  $X=(X_t)_{t\in \Z}\in {\cal
 M}_\Z(f_\theta,M_\theta)$ is assured (see \cite{Bardet2009}). The subset $\Theta(r)$ is defined as a reunion to consider
 accurately general causal models and ARCH-type models simultaneously.\\
 The following assumptions are needed to study QLME property.\\
 \noindent {\bf Assumption D$(\Theta)$:} $\exists\underline{h}>0$ such that
$\displaystyle \inf_{\theta \in
 \Theta}(|h_\theta(x)|)\geq \underline{h}$ for all $x\in \R^{\N}.$ \\
~\\
 {\bf Assumption Id($\Theta$):} For all  $\theta,\theta'\in \Theta^2$,
$$ \Big( f_{\theta}(X_0,X_{-1},\cdots)=f_{\theta'}(X_0,X_{-1},\cdots)~\mbox{and}~h_{\theta}(X_0,X_{-1},\cdots)=h_{\theta'}(X_0,X_{-1},\cdots) \ \text{a.s.}\Big) \Rightarrow \ \theta = \theta'.$$
{\bf Assumption Var($\Theta$):} For all $\theta  \in \Theta $, one
 of the families $ \big( \dfrac{\partial f_{\theta}}{\partial \theta^{i}}(X_0,X_{-1},\cdots) \big)_{1\leq i \leq d} \quad
  \mbox{or} \quad \big( \dfrac{\partial h_{\theta}}{\partial \theta^{i}}(X_0,X_{-1},\cdots) \big)_{1\leq i \leq d}  \quad $
 is a.s. linearly independent.\\

     As in \cite{Bardet2009}, we will make the convention that if \textbf{A}$_i(M_\theta,\Theta)$ holds then $\alpha_\ell^{(i)}( h_\theta,\Theta) = 0$ and if
     \textbf{A}$_i(h_\theta,\Theta)$ holds then $\alpha_\ell^{(i)}( M_\theta,\Theta) = 0.$ Denote :\\

 \noindent \textbf{ Assumption K}($f_\theta,M_\theta,\Theta$) : for i= 0, 1, 2, \textbf{A}$_i(f_\theta,\Theta)$ and \textbf{A}$_i(M_\theta,\Theta)$
 (or \textbf{A}$_i(h_\theta,\Theta)$) hold and there exists $l>2$ such that $ \alpha_j^{(i)}(f_{\theta},\Theta) + \alpha_j^{(i)}(M_{\theta},\Theta) + \alpha_j^{(i)}(h_{\theta},\Theta) =
      \mathcal{O}(j^{-l})$,  for i= 0, 1.\\
 Throughout the sequel, we will assume that the functions $ \theta \mapsto M_\theta$ and $\theta \mapsto f_\theta$ are twice continuously differentiable on
 $\Theta$.

 \subsection{Examples}

 \begin{enumerate}
    \item AR($\infty$) models.\\
          Consider the AR($\infty$) process defined by :
          $$ X_t = \sum_{k\geq 1} \phi_k(\theta^*_0)X_{t-k} + \xi_t  ~ , ~ t \in \Z  $$
 with $ \theta^*_0 \in \Theta$, where $\Theta$ is a compact subset of $\R^d$ such that $ \sum_{k\geq 1} \| \phi_k(\theta) \|_{\Theta} < 1$.
  The process belongs to the class ${\cal M}_{\Z}(M_{\theta^*_0},f_{\theta^*_0})$ where
  $f_\theta(x_1,\cdots)= \sum_{k\geq 1} \phi_k(\theta)x_k $ and $ M_\theta \equiv 1 ~ \text{for all}~ \theta \in \Theta$.
  Then Assumptions D($\Theta$) and $A_0(f_\theta,\Theta)$ hold  with $\underline{h}=1$ and $\alpha^{(0)}_k(f_\theta,\Theta) = \|\phi_k(\theta)\|_\Theta$.
 If there exists $\ell >2$ and $\phi_k$ twice differentiable such as  $ \|\phi_k(\theta)\|_\Theta =\|\phi_k'(\theta)\|_\Theta = \|\phi_k''(\theta)\|_\Theta  = O(k^{-\ell})$, then
   Assumptions \textbf{K}($f_\theta,M_\theta,\Theta$) holds. Moreover, if $\xi_0$ is a nondegenerate random variable,
   Id($\Theta$) and   Var($\Theta$) hold. For any $r\geq 1$ such that $ \E|\xi_0|^r<\infty$, $\Theta(r)=\Theta$.

    \item GARCH(p,q) models.\\
          Consider the GARCH(p,q) process defined by :
 $$ X_t= \sigma_t \, \xi_t \ ,\ \sigma_t^2 = \alpha_{0}^*+ \sum^{q}_{k=1} \alpha_{k}^*X^2_{t-k} + \sum^{p}_{k=1} \beta_{k}^*\sigma^2_{t-k}  \ ,  \ \  t\in \Z $$
 with $\E(\xi_0^2)=1$ and $ \theta^*_0 := (\alpha^*_0,\cdots,\alpha^*_q,\beta^*_1,\cdots,\beta^*_p ) \in \Theta $ where $\Theta$ is a
 compact subset of $ ]0 , \infty[\times [0 , \infty [^{p+q}$   such that
 $\sum_{k=1}^{q}\alpha_k +  \sum_{k=1}^{p}\beta_k < 1 $ for all $\theta \in \Theta$.
 %
 Then there exists (see Bollerslev \cite{Bollerslev1986} or Nelson and Cao \cite{Nelson1992}) a nonnegative sequence $(\psi_k(\theta^*_0))_{k\geq 0}$
 such that $\sigma_t^2 = \psi_0(\theta^*_0) + \sum^{}_{k\geq 1} \psi_k(\theta^*_0)X^2_{t-k}$  with  $ \psi_0(\theta^*_0) = \alpha^*_0/(1-\sum^{p}_{k=1} \beta_{k}^*)$.
 This process belongs to a class ${\cal M}_{\Z}(M_{\theta^*_0},f_{\theta^*_0})$ where
  $M_\theta(x_1,\cdots)= \sqrt{\psi_0(\theta) + \sum_{k\geq 1} \psi_k(\theta)x_k} $ and $ f_\theta \equiv 0 ~ \forall \theta \in \Theta$.
   Assumptions D($\Theta$)  holds with $\underline{h}= \underset{\theta\in \Theta} {\mbox{inf}}(\alpha_0)$.
  If there exists $0<\rho_0<1$ such that for any $\theta \in \Theta, ~ \sum_{k=1}^{q}\alpha_k +  \sum_{k=1}^{p}\beta_k \leq \rho_0$ then
 the sequences $( \|\psi_k(\theta)\|_\Theta)_{k\geq 1}$, $( \|\psi'_k(\theta)\|_\Theta)_{k\geq 1}$ and $( \|\psi''_k(\theta)\|_\Theta)_{k\geq 1}$
  decay exponentially fast (see Berkes {\it et al.} \cite{Berkes2003}), thus Assumption \textbf{K}($f_\theta,M_\theta,\Theta$) holds.
  Moreover, if $\xi^2_0$ is a nondegenerate random variable,
   Id($\Theta$) and   Var($\Theta$) hold.  For $r\geq 2$ denote
   $$  \Theta(r)= \big\{  \theta \in \Theta ~ ; ~  (\E|\xi_0|^r)^{2/r} \sum_{k=1}^{q}\alpha_k +  \sum_{k=1}^{p}\beta_k < 1 \big\}  .$$

 \end{enumerate}

\subsection{Test statistics}

       Assume that a trajectory $ (X_1, \cdots, X_n)$ is
       observed. It is clear that if $ (X_1, \cdots, X_n) \in \mathcal{M}_{ \{1,\cdots,n\} }(M_\theta,f_\theta)$, then
        for $ T\subset \{1,\cdots,n \} $,  the conditional quasi-(log)likelihood computed on $T$ is given by :
       $$ L_n(T,\theta) := -\dfrac{1}{2} \sum\limits_{t \in T} q_t(\theta)   ~~ \text{with}  ~~ q_t(\theta)= \dfrac{(X_t-f^t_{\theta})^2}{h_\theta^t} + \log(h_\theta^t) $$
       where $ f^t_{\theta}=f_{\theta}\big(X_{t-1},X_{t-2}\ldots\big) $, $ M^t_{\theta}=M_{\theta}\big(X_{t-1},X_{t-2}\ldots\big) $ and
   $ h^t_{\theta}= {M^{t}_{\theta}}^2 $. Therefore, we approximate
   the conditional log-likelihood with :

   $$ \widehat L_n(T,\theta):=-\frac{1}{2}\sum\limits_{t \in T}\widehat q_t(\theta)\quad
 \textrm{where}\;\;\;\widehat{q}_t(\theta):=\frac{\big(X_t-\widehat{f}^t_{\theta}\big)^2}{\widehat{h}{^t_{\theta}}} +\log\big(\widehat{h}{^t_{\theta}}\big) $$
 with
 $ \widehat{f}^t_{\theta}=f_{\theta}\big(X_{t-1},\ldots,X_{1},0,0,\cdots\big)$, $ \widehat{M}^t_{\theta}=
 M_{\theta}\big(X_{t-1},\ldots,X_{1},0,0,\cdots\big) $ and $\widehat{h}^t_{\theta}=( \widehat{M}^{t}_{\theta} )^2$.\\
   For $ T\subset \{1,\cdots,n \} $, define the estimator $\widehat{\theta}_n(T) :=  \underset{\theta\in \Theta}
   {\mbox{argmax}}(\widehat{L}_n(T,\theta))$.
   Moreover, for $ 1 \leq k \leq n $, denote $ T_k= \{1,\cdots,k \}$ and $ \overline{T}_k= \{k+1,\cdots,n \}$.\\

  \noindent   Now, define\\
   $\widehat{G}_n(T) := \dfrac{1}{Card(T)} \sum \limits_{t \in T} \Big( \dfrac{\partial \widehat{q}_{t}( \widehat{\theta}_n(T))}{\partial \theta} \Big)
                       \Big( \dfrac{\partial \widehat{q}_{t}( \widehat{\theta}_n(T))}{\partial \theta} \Big)' $
   and
   $ \widehat{F}_n(T):= - \dfrac{2}{Card(T)} \Big(  \dfrac{\partial^2 \widehat{L}_{n} (T,\widehat{\theta}_n(T))}{\partial \theta \partial \theta'}
   \Big)$.\\

   \noindent  For $k=1,\cdots,n-1$, denote :
   $$ \widehat{\Sigma}_{n,k} := \dfrac{k}{n}\widehat{F}_n(T_k)\widehat{G}_n(T_k)^{-1}\widehat{F}_n(T_k)\text{\large{\1}}_{\det(\widehat{G}_n(T_k))\neq 0 }
      +  \dfrac{n-k}{n}\widehat{F}_n(\overline{T}_k)\widehat{G}_n(\overline{T}_k)^{-1}\widehat{F}_n(\overline{T}_k)\text{\large{\1}}_{\det(\widehat{G}_n(\overline{T}_k))\neq 0 } .$$
 Let $(v_n)_{n\in \N}$  be a sequence satisfying  $ v_n \to \infty $ and $v_n/{n} \to 0$ (as $n\to
 \infty$). Denote $\Pi_n= [v_n, n-v_n]\cap \N  $ and define the  statistics:\\

 $  \widehat{Q}^{(1)}_n:= \underset{ k \in \Pi_n } {\mbox{max}}\widehat{Q}^{(1)}_{n,k} ~~ ~
 \text{where} ~~  \widehat{Q}^{(1)}_{n,k}:= \dfrac{k^2}{n} \big( \widehat{\theta}_n(T_k)-  \widehat{\theta}_n(T_n) \big)' \widehat{\Sigma}_{n,k} \big( \widehat{\theta}_n(T_k)-  \widehat{\theta}_n(T_n) \big);$\\

   $  \widehat{Q}^{(2)}_n:= \underset{ k \in \Pi_n } {\mbox{max}}\widehat{Q}^{(2)}_{n,k} ~~ ~
 \text{where} ~~  \widehat{Q}^{(2)}_{n,k}:= \dfrac{(n-k)^2}{n} \big( \widehat{\theta}_n(\overline{T}_k)-  \widehat{\theta}_n(T_n) \big)' \widehat{\Sigma}_{n,k} \big( \widehat{\theta}_n(\overline{T}_k)-  \widehat{\theta}_n(T_n) \big);$\\

   $ \widehat{Q}_n:=\text{max} \big(\widehat{Q}^{(1)}_n , \widehat{Q}^{(2)}_n  \big)$ which is the test statistic.\\

\begin{rem}


    Note that,  the choice of $(v_n)$ is crucial in practice.
    We evaluated the procedure with $ v_n=[\log n], ~ [(\log n)^2], ~ [(\log n)^3] $ and recommend to use $ v_n= [(\log n)^2]$
   for linear model and $ v_n= [(\log n)^{5/2}]$ for GARCH-type model.

\end{rem}

   Lee and Song \cite{Lee2008} constructed a test for detecting changes in parameters of ARMA-GARCH models.
   Their test statistic uses a matrix $ \widehat{\Sigma}_n $ which depends on the estimator $\widehat{\theta}_n(T_n)$.
   Under the null hypothesis (the parameter $\theta_0$ does not change), $ \widehat{\Sigma}_n $ is a consistent estimator of $FG^{-1}F$ where
    $G= \E [(\partial q_0(\theta_0)/\partial \theta) (\partial q_0(\theta_0)/\partial \theta)'  ]$  and $F=\E [ \partial ^2 q_0(\theta_0)/\partial \theta \partial \theta' ] $.
    Under the alternative, the model depends on several parameters and  $\widehat{\theta}_n(T_n)$  may not be a consistent estimator of one of them.
    Therefore, the convergence of the matrix $ \widehat{\Sigma}_n $ is not assured.
    Thus, the asymptotic behavior of the test statistic may be very difficult to study.
    To solve this problem, we introduce the family of matrices  $ \{ \widehat{\Sigma}_{n,k}, ~ k\in \Pi_n  \}  $.
    It is easy to see that under the null hypothesis, any sequence $ (\widehat{\Sigma}_{n,k_n})_{n > 1 ,k_n \in \Pi_n}$ is consistent.
    We show (see proof of Theorem \ref{theo2}) that under the local alternative where there is one change in the model,
    there exists a sequence $ (\widehat{\Sigma}_{n,k^*_n})_{n > 1 ,k^*_n \in \Pi_n}$ which converges.

 \section{ Asymptotic results}
   \subsection{ Asymptotic under the null hypothesis}

\begin{Theo}\label{theo1}
 Assume \textbf{D}$(\Theta)$, \textbf{Id}($\Theta$),
 \textbf{Var} and \textbf{K}($f_\theta,M_\theta,\Theta$). Under the null
 hypothesis $H_0$, if $ \theta_0 \in \overset{\circ}{\Theta}(4)$, then for
  $j=1,2$,
 $$  \widehat{Q}^{(j)}_n \limiteloin \underset{ 0 \leq \tau \leq 1 } {\mbox{sup}} \| W_d(\tau) \|^2  $$
where $W_d$ is a d-dimensional Brownian bridge.

\end{Theo}

 For any $\alpha \in (0,1)$, let $C_\alpha$ denote the $(1-\alpha/2)$-quantile of the distribution of
  $ \underset{ 0 \leq \tau \leq 1 } {\mbox{sup}} \| W_d(\tau)\|^2$. Then, the following corollary is a direct application of Theorem \ref{theo1}.

 \begin{cor}\label{cor1}
  Under assumptions of Theorem \ref{theo1} :
  $$ \forall \alpha \in (0,1) ~~   \underset{ n \rightarrow \infty  } {\limsup }~ P \big( \widehat{Q}_n > C_\alpha \big) \leq \alpha.   $$
 \end{cor}

\begin{rem}\label{rem1}
     Quantile values of the distribution of $ \underset{ 0 \leq \tau \leq 1 } {\mbox{sup}} \| W_d(\tau)\|^2$
   are known (see for instance Kieffer \cite{Kiefer1959} for $d \in \{1,\cdots,5 \}$ or Lee {\it et al.} \cite{Lee2003} for $d \in \{1,\cdots,10 \}$).
\end{rem}


 Theorem \ref{theo1} and Corollary \ref{cor1} imply that a large value of $\widehat{Q}_n$  means there is a
 change in the model. At a nominal level $\alpha$, the critical region of the test is  $(\widehat{Q}_n > C_\alpha )$.\\

    Figure 1 is an illustration of the test procedure for AR(1) process. At a level $ \alpha = 0.05$, for $d=1$,  $ C_\alpha \simeq 2.20$. Figure 1 a-) and b-) show that,
    the values of $\widehat{Q}^{(1)}_{n,k}$ and $\widehat{Q}^{(2)}_{n,k}$ are all below the red line which represents the limit of the critical region.
    Figure 1 c-) and d-) show that  $\widehat{Q}^{(1)}_{n,k}$ and $\widehat{Q}^{(2)}_{n,k}$ are larger and  increase around the point where the change occurs.\\
    \indent As it can be observed on the Figure 1 and Figure 2 , the statistics $\widehat{Q}^{(1)}_{n}$ and $\widehat{Q}^{(1)}_{n}$ are not clearly equal.
    Figure 2 shows the typical example for $ARCH(1)$ with one change where $ \widehat{Q}^{(1)}_n < C_\alpha$ and $ \widehat{Q}^{(2)}_n > C_\alpha$.
    In general, we don't know if under the alternative hypothesis each of statistics $\widehat{Q}^{(1)}_n$ and $  \widehat{Q}^{(2)}_n$ take large values.
    But, under the local alternative of one change, the maximum of these two statistics diverges to infinity (see Theorem \ref{theo2}).
    This is the reason why we define the critical region  as $\{ \text{max}(\widehat{Q}^{(1)}_n , \widehat{Q}^{(2)}_n )> C_\alpha \}  $.

\subsection{ The asymptotic under a local alternative}
  In this subsection, we consider a local alternative that there is one
  change in the model. More precisely, define
\begin{multline*}
   H_1^{(loc)}:  ~ \text{there exist} ~ ~ \tau^* \in (0 , 1) ~ \text{and} ~  \theta^*_1, \theta^*_2 \in \Theta ~ \text{with}~ \theta^*_1\neq\theta^*_2
   ~ \text{such that} ~ X_1, \cdots X_{[n\tau^*]} \in \mathcal{M}_{T_{[n\tau^*]}}(M_{\theta^*_1},f_{\theta^*_1})\\
    ~ \text{and} ~ X_{[n\tau^*]+1},\cdots,X_n \in\mathcal{M}_{\overline{T}_{[n\tau^*]}}(M_{\theta^*_2},f_{\theta^*_2}).\hspace{8cm}
 \end{multline*}
 \begin{Theo}\label{theo2}
 Assume \textbf{D}$(\Theta)$, \textbf{Id}($\Theta$),
 \textbf{Var} and \textbf{K}($f_\theta,M_\theta,\Theta$). Under $H_1^{(loc)}$, if ~$ \theta^*_1, \theta^*_2 \in \overset{\circ}{\Theta}(4)$, then
 $$  \widehat{Q}_n \limitepsn \infty  .$$
\end{Theo}
\begin{rem}
     1-) Theorem \ref{theo2} shows that the test with the local alternative $ H_1^{(loc)}$  is consistent in power.\\
    2-) This procedure can also be used  to test multiple changes using  ICSS type algorithm developed by Inclán and Tiao \cite{Tiao1994}.
\end{rem}
%

\[
 \epsfxsize 16cm \epsfysize 12cm \epsfbox{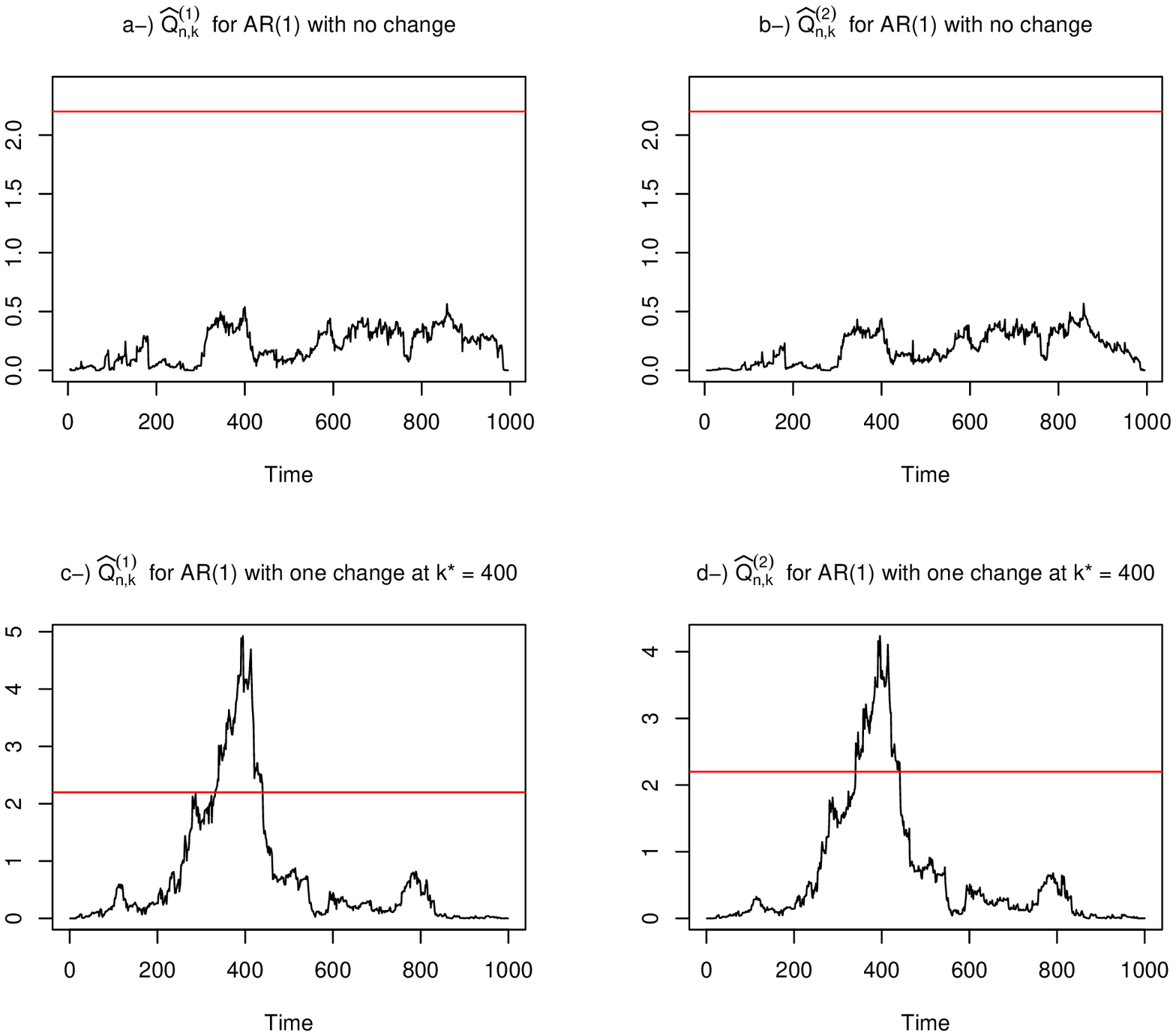}
\vspace{-.5cm}
\]
 {\bf Figure 1 :}   Typical realization of the statistics  $\widehat{Q}^{(1)}_{n,k}$ and $\widehat{Q}^{(2)}_{n,k}$ for 1000 sample of AR(1) with $v_n= [\log n]$.
  a-) and b-) are the
 case of AR(1) where parameter $\phi_1=0.3$ remains constant. c-) and d-)
 are the case of AR(1) with parameter $\phi_1=0.3$  changing to $ 0.5 $ at $k^* = 400.$  \\
~\\

%
%

\vspace{-1.5cm}
\[
 \epsfxsize 15cm  \epsfysize 8cm \epsfbox{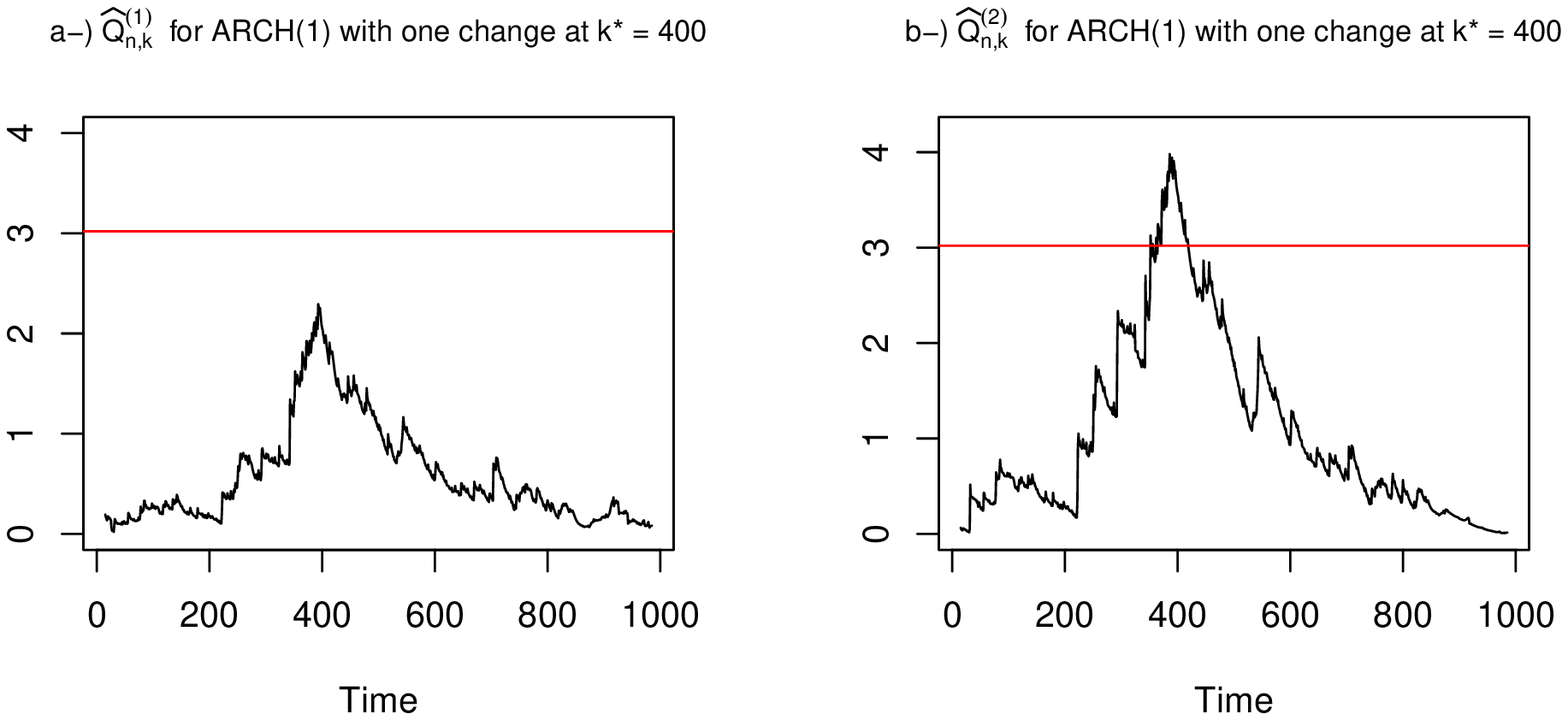}
 \vspace{-0.5cm}
\]
 {\bf Figure 2 :}  Typical realization of the statistics  $\widehat{Q}^{(1)}_{n,k}$ and $\widehat{Q}^{(2)}_{n,k}$ for 1000 sample of ARCH(1)  with $v_n= [(\log n)^2]$.
 a-) and b-)  are the case of ARCH(1) with parameter  $\theta_1=(1,0.3)$  changing to $(1,0.15)$ at $k^* = 400$.\\
~\\

\section{ Some simulations results}

 In this section, we evaluate the performance of the procedure through empirical study. We compare our results with those obtained
 by Kouamo {\it et al.} \cite{Kouamo2010}, Lee and Na \cite{Lee2005} and the results obtained from the residual CUSUM test by using the
 statistics defined by Kulperger and Yu \cite{Kulperger2005}. For a sample size $n$, $ \widehat{Q}_n $ is computed  with $v_n= [(\ln n)^2]$ for AR model
 and $v_n= [(\ln n)^{5/2}]$ for GARCH model  and is compared to the critical value of the test.
 In the following models, $(\xi_t)_{t \in \Z}$ are iid standard Gaussian random variables.

 \subsection{ Test for parameter change in AR(p) models}
  Let us consider a AR(p) process : $X_t= \sum\limits_{k=1}^{p}\phi^*_k X_{t-k} +
  \xi_t$ with $p \in \N^* $. The true parameter of
  the model is denoted by $ \theta^*_0 = (\phi^*_1, \cdots, \phi^*_p) \in \Theta$ where
  $ \Theta = \{  \theta = ( \phi_1, \cdots,  \phi_p) \in \R^p ~ ~ /  \sum\limits_{i=1}^{p}| \phi_i|<1   \}   $.
  Since $ M_\theta \equiv 1$, $ \Theta (r) = \Theta $ for any $r\geq 1$. Assume $(X_1, \cdots, X_n)$ is observed,
  we have for any $\theta \in \Theta$, $ \widehat{q}_t(\theta)= \big(X_t- \sum\limits_{k=1}^{p}\phi_k X_{t-k} \big)^2$,
   $  \dfrac{\partial \widehat{q}_{t}( \theta)}{\partial \theta}  =
  -2 \big(X_t- \sum\limits_{k=1}^{p}\phi_k X_{t-k} \big)\cdot (X_{t-1}, X_{t-2},\cdots,X_{t-p})$
   and for $  1 \leq i,j\leq n$ ~
  $ \dfrac{\partial^2 \widehat{q}_{t}( \theta)}{\partial \phi_i \partial \phi_j} = 2X_{t-i} X_{t-j}  .$\\

  We consider a  AR(1) process with one parameter.  At level $\alpha = 0.05$,  the critical value is $C_{\alpha} \simeq 2.20 $.
   For $n=1024, 2048, 4096$ ;  we generate a sample $(X_1,\cdots,X_n)$ in the following situations :
   (i) there is no change and the parameter of the model is $ \theta_0= 0.9 $ and
    (ii) there is one change and the parameter $ \theta_0= 0.9$ changes to $ \theta_1=0.5$ at $n/2$. The following table indicates the
    proportion of the number of rejections of the null hypothesis out of 100 repetitions.

   \begin{table}[htbp]
    \centering
    \begin{tabular}{|c|c|c|c|}
      \hline
        &  $n=1024$ &  $n=2048$ & $n=4096$ \\
      \hline
      \hline
      &     &      &   \\
   Empirical levels & 0.080  ~ (0.134 ; 0.092)  &  0.070  ~ (0.100 ; 0.062) &  0.050  ~ (0.082 ; 0.040) \\
     &        &      &      \\
     Empirical powers when $ \theta_1=0.5 $  &  0.980  ~ (0.590 ; 0.530) &  0.990 ~ (0.720 ; 0.680) & 0.990  ~ (0.810 ; 0.790) \\
      \hline
    \end{tabular}
 \caption{{\footnotesize  Empirical levels and powers at nominal level 0.05 of the test for parameter change in AR(1) model.
 The empirical levels are computed when $ \theta_0= 0.9$ ; the empirical powers are computed when $ \theta_0$ changes to $ \theta_1=0.5$ at $n/2$.
   Figures in brackets are the results obtained by Kouamo {\it et al.} \cite{Kouamo2010} at the scale J=4 with KSM and CVM statistic in wavelet domain.}}
        \label{tab1}
\end{table}
    Table \ref{tab1}  shows that the empirical level of the test decreases as n increases and equals to 0.05 when n = 4096.
    These levels are close to those obtained by Kouamo {\it et al.} with CVM (Cramér-Von Mises) test statistic.
    The results obtained with our test statistic $\widehat{Q}_n$ are clearly more accurate.

  \subsection{ Test for parameter change in GARCH(1,1) models}
  Consider the GARCH(1,1) model defined by:
   $$ \forall t \in \Z, ~ ~   X_t= \sigma_t \xi_t ~~ \text{ with} ~~ \sigma^2_t= \alpha^*_0 + \alpha^*_1 X^2_{t-1} + \beta^*_1\sigma^2_{t-1}      $$
  with $ \theta^*_0=(\alpha^*_0,\alpha^*_1,\beta^*_1) \in \Theta \subset ]0,\infty[\times[0,\infty[^2$ and satisfying $ \alpha^*_1 + \beta^*_1 < 1 $.
  The ARCH($\infty$) representation is
    $  \sigma^2_t = \alpha^*_0/(1-\beta^*_1) + \alpha^*_1 \sum\limits_{ k\geq 1}^{} (\beta^*_1)^{k-1} X^2_{t-k}.  $
 %
  For any $ \theta \in \Theta  $  and $ t = 2, \cdots, n$ , we have
   $$ \widehat{h}^t_{\theta} = \alpha_0/(1-\beta_1) + \alpha_1X^2_{t-1}  + \alpha_1 \sum\limits_{ k= 2}^{t} \beta_1^{k-1} X^2_{t-k}  ~~ \text{and} ~~
      \widehat{q}_t(\theta) =  X^2_t /~ \widehat{h}^t_{\theta}  +  \log(\widehat{h}^t_{\theta})  .$$
    Therefore, it follows that
   $ \dfrac{\partial \widehat{q}_t(\theta)}{ \partial \theta}  =  \dfrac{1}{\widehat{h}^t_{\theta} } \Big( 1 -\dfrac{X^2_t}{\widehat{h}^t_{\theta} } \Big)
      \Big( \dfrac{\partial \widehat{h}^t_{\theta}}{ \partial \alpha_0}  , \dfrac{\partial \widehat{h}^t_{\theta}}{ \partial \alpha_1}   ,
        \dfrac{\partial \widehat{h}^t_{\theta}}{ \partial \beta_1} \Big)   $
   with  $ \partial \widehat{h}^t_{\theta} /  \partial \alpha_1 = X^2_{t-1} + \sum\limits_{ k= 2}^{t} \beta_1^{k-1} X^2_{t-k} $
   $  \partial \widehat{h}^t_{\theta} /  \partial \alpha_0 = 1/(1-\beta_1) $,
   and ~ $ \partial \widehat{h}^t_{\theta} /  \partial \beta_1 = \alpha_0/ (1-\beta_1)^2 + \alpha_1 X^2_{t-2} +  \alpha_1 \sum\limits_{ k= 3}^{t} (k-1)\beta_1^{k-2} X^2_{t-k} $.

    \noindent Let $ \theta = (\alpha_0,\alpha_1,\beta_1) = (\theta_1,\theta_2,\theta_3) \in \Theta  $, for $1 \leq i,j \leq 3$, we have
    $$ \dfrac{\partial^2 \widehat{q}_t(\theta)}{ \partial \theta_i  \partial \theta_j} =  \dfrac{1}{(\widehat{h}^t_{\theta})^2 }
  \Big( \dfrac{2X^2_t}{\widehat{h}^t_{\theta} } -1 \Big) \dfrac{\partial \widehat{h}^t_{\theta}}{ \partial \theta_i} \dfrac{\partial \widehat{h}^t_{\theta}}{ \partial \theta_j}
  +\dfrac{1}{\widehat{h}^t_{\theta} } \Big( 1 -\dfrac{X^2_t}{\widehat{h}^t_{\theta} } \Big) \dfrac{\partial^2 \widehat{h}^t_{\theta}}{\partial \theta_i \partial \theta_j } $$
   with ~ $ \partial^2 \widehat{h}^t_{\theta} / \partial \alpha_0^2  = 0$,  $\partial^2 \widehat{h}^t_{\theta} / \partial \alpha_0  \partial \alpha_1 =0$,
     $\partial^2 \widehat{h}^t_{\theta} / \partial \alpha_1^2  = 0 $,
     $ \partial^2 \widehat{h}^t_{\theta} / \partial \alpha_1  \partial \beta_1  = X^2_{t-2} + \sum\limits_{ k= 3}^{t} (k-1)\beta_1^{k-2} X^2_{t-k} $,
    $ \partial^2 \widehat{h}^t_{\theta} / \partial \alpha_0  \partial \beta_1  = 1/(1-\beta_1)^2 $
     and ~ $ \partial \widehat{h}^t_{\theta} /  \partial \beta^2_1 = 2 \alpha_0/ (1-\beta_1)^3 +2\alpha_1X^2_{t-3} +  \alpha_1 \sum\limits_{ k= 4}^{t} (k-1)(k-2)\beta_1^{k-3} X^2_{t-k} $.\\                          ~

  \begin{enumerate}
    \item  Case of ARCH(1). Assume $\beta_1=0$ and $\theta = (\alpha_0,\alpha_1)$. At level $\alpha = 0.05$,  the critical value is $C_{\alpha}  \simeq 3.02 $.
    For $n=500, 800, 1000$ ;  we generate a sample $(X_1,\cdots,X_n)$ in the following situations  :
   (i) there is no change, the parameter of the model is $ \theta_0= (1, 0.3)$ and
    (ii) there is one change, the parameter $ \theta_0= (1, 0.3)$
    changes to two different values of $ \theta_1=(0.5,0.3)$ and $\theta_1=(0.5,0.6)$ at $n/2$. The following table indicates the proportion of the number of rejections
   of the null hypothesis out of 500 repetitions.

   \begin{table}[htbp]
    \centering
    \begin{tabular}{|c|c|c|c|}
      \hline
        &  $n=500$ &  $n=800$ & $n=1000$ \\
      \hline
      \hline
      &    &    &    \\
   Empirical levels & 0.068  ~ (0.088)  &  0.048  ~ (0.080) &  0.036  ~ (0.074) \\
       &     &    &     \\
     Empirical powers when $ \theta_1=(0.5, 0.3)$  &  0.948  ~ (0.922) & 0.972  ~ (0.987) & 0.998  ~ (0.995) \\
        &     &     &    \\
     Empirical powers when $ \theta_1=(0.5, 0.6)$  &  0.876  ~ (0.498) &  0.976  ~ (0.589) &  0.984 ~ (0.606)  \\
      \hline
    \end{tabular}
 \caption{{\footnotesize  Empirical levels and powers at nominal level 0.05  of the test for parameter change in ARCH(1) model.
 The empirical levels are computed when $ \theta_0= (1, 0.3)$ ; the empirical powers are computed when $ \theta_0$ changes to $ \theta_1$ at $n/2$.
   Figures in brackets are the results obtained by Lee and Na \cite{Lee2005}.}}
        \label{tab2}
\end{table}

    \item Case of GARCH(1,1). Now $\theta = (\alpha_0,\alpha_1,\beta_1)$.
    At level $\alpha = 0.05$, the critical value is $C_{\alpha} \simeq 3.47 $.
    For $n=500, 1000, 1500$ ;  we generate a sample $(X_1,\cdots,X_n)$ in the following situations :
   (i) there is no change, the parameter of the model is $ \theta_0= (1, 0.4, 0.1)$ and
    (ii) there is one change, the parameter $ \theta_0= (1, 0.4, 0.1)$ changes two different values of $ \theta_1=(0.7,0.4,0.1)$ and $\theta_1=(1,0.4,0.3)$ at $n/2$.
    The following table indicates the
    proportion of the number of rejections of the null hypothesis out of 500 repetitions.

  \end{enumerate}

  \begin{table}[htbp]
    \centering
    \begin{tabular}{|c|c|c|c|}
      \hline
        &  $n=500$ &  $n=1000$ & $n=1500$ \\
      \hline
      \hline
      &     &    &    \\
   Empirical levels & 0.100  ~ (0.030)  &  0.078  ~ (0.032) &  0.052  ~ (0.042) \\
   &       &    &   \\
     Empirical powers when $ \theta_1=(0.7, 0.4, 0.1)$  &  0.498  ~ (0.334) & 0.752 ~ (0.658) & 0.934  ~ (0.848) \\
     &     &    &    \\
    Empirical powers when $ \theta_1=(1, 0.4, 0.3)$  &  0.654  ~ (0.404) &  0.968  ~ (0.772) &  0.976 ~ (0.922)  \\
      \hline
    \end{tabular}
 \caption{{\footnotesize  Empirical levels and powers at nominal level 0.05  of test the for parameter change in GARCH(1,1) model.
 The empirical levels are computed when $ \theta_0= (1, 0.4, 0.1)$ ; the empirical powers are computed when $ \theta_0$ changes to $ \theta_1$ at $n/2$.
   Figures in brackets are the results of the residual CUSUM test using  $CUSUM^{(2)}$ statistic defined by Kulperger and Yu \cite{Kulperger2005}. }}
        \label{tab3}
\end{table}

    Table \ref{tab2}  and Table \ref{tab3} show that the empirical level of the test decreases and the empirical power increases
     as $n$ increases. For ARCH model, we can see that the empirical level is less than 0.05 when $n=800$.
     It is not very surprising because the asymptotic size of the test is less than $\alpha=0.05$.
      This is not the case for GARCH model.
    It is explained by the fact that the application of the procedure  to GARCH model requires ARCH($\infty$) representation.
    Thus, the information contained in all the past of the process is not used because it is not observed.
    In Table \ref{tab2}, figures in brackets are the results obtained by Lee and Na \cite{Lee2005} using the CUSUM test based
    on conditional least-squares estimator.
     In Table \ref{tab3}, figures in brackets are the results of the residual CUSUM test that we obtained by using  $CUSUM^{(2)}$ statistic studied by
      Kulperger and Yu \cite{Kulperger2005}. Once again, our test statistic $\widehat{Q}_n$  provides best results.

 \section{ Proofs of the main results}

  Let $ (\psi_n)_n  $ and $ (r_n)_n  $ be sequences of random variables. Throughout this section, we use the notation
  $ \psi_n = o_P(r_n)  $ to mean :  for all $  \varepsilon  > 0, ~ P( |\psi_n| \geq \varepsilon|r_n| ) \rightarrow 0 $ as $n \rightarrow  \infty$.
  Write $ \psi_n = O_P(r_n)  $ to mean :  for all  $  \varepsilon > 0 $,  there exists $C>0$  such that  $ ~ P( |\psi_n| \geq C |r_n| )\leq \varepsilon  $
   for $n$ large enough.

\subsection{Some preliminary results}
  First, let us prove   useful technical lemmas.

  Under the null hypothesis $H_0$  the observations $ (X_1, \cdots, X_n)$ belong in the
       class $\mathcal{M}_{ \{1,\cdots,n \} }(M_{\theta_0},f_{\theta_0})$, define the matrix
  $G:= \E \Big [ \dfrac{\partial q_0(\theta_0)}{\partial \theta} \dfrac{\partial q_0(\theta_0)}{\partial \theta} ' \Big]$
  ( where $'$ denotes the transpose) and
  $ F:=  \E \Big [ \dfrac{\partial^2 q_0(\theta_0)}{\partial \theta \partial \theta'}  \Big] .$
  Under assumption \textbf{Var}, F is a non-singular matrix (see \cite{Bardet2009}).

  \begin{lem} \label{lem_inv_G}
  Assume the functions $ \theta \mapsto M_\theta $  and $ \theta \mapsto f_\theta$
  are 2-times continuously differentiable on $\Theta$. Under the null hypothesis \textbf{D}$(\Theta)$ and \textbf{Var}, G  is a symmetric,
  positive definite matrix.
 \end{lem}
\vspace{-0.5cm}
\begin{dem} \label{lem_inv_G}
 It is clear that G is symmetric.
  Moreover, for $1\leq i \leq d$, we have :\\
  $\dfrac{\partial q_{0}(\theta_0)}{\partial \theta_i} = -2\dfrac{ \xi_0}{\sqrt{h^0_{\theta_0}}} \dfrac{\partial
  f^0_{\theta_0}}{\partial \theta_i}- \dfrac{ \xi^2_0}{h^0_{\theta_0}} \dfrac{\partial
  h^0_{\theta_0}}{\partial \theta_i}+ \dfrac{1}{h^0_{\theta_0}} \dfrac{\partial h^0_{\theta_0}}{\partial
  \theta_i}.$
  Thus, using independence of $\xi_0$ and $X_{-1},~ X_{-2},\cdots $   we obtain :
  \begin{equation}\label{mat_G_sum}
  \E \Big[ \dfrac{\partial q_0(\theta_0)}{\partial \theta}' \dfrac{\partial q_0(\theta_0)}{\partial \theta}  \Big]
     = 4 \E \Big [ \dfrac{ 1}{h^0_{\theta_0}} \dfrac{\partial f^0_{\theta_0}}{\partial \theta}' \dfrac{\partial f^0_{\theta_0}}{\partial \theta}  \Big]
     + \E \Big( (\xi^2_0-1)^2 \Big ) \E \Big [ \dfrac{ 1}{(h^0_{\theta_0})^2} \dfrac{\partial h^0_{\theta_0}}{\partial \theta}' \dfrac{\partial h^0_{\theta_0}}{\partial \theta}  \Big]
     .
      \end{equation}
     Since $\E \xi^2_0 = 1$, it is easy to see that $\E \Big( (\xi^2_0-1)^2 \Big )
     >0.$

   Under \textbf{Var}, one of the two matrix of the right-hand side
   of relation (\ref{mat_G_sum}) is positive definite and the other
  is semi-positive definite. Thus, G is positive definite.
 \end{dem}

Now, recall that $ F:=  \E \Big [ \dfrac{\partial^2
q_0(\theta_0)}{\partial \theta \partial \theta'}  \Big]$.
  Let $ T\subset \{1,\cdots,n \}.$  For any $ \theta \in \Theta $ and  $ i=1,\cdots,d $, by Taylor expansion of $ \partial L_n(T,\theta_0)/\partial \theta_i $,
  there exist $ \overline{\theta}_{n,i} \in [\theta_0,  \theta] $
 such that:
 \begin{equation}\label{Taylor_simple}
   \dfrac{\partial L_n(T, \theta)}{\partial \theta_i}  =  \dfrac{\partial L_n(T, \theta_0)}{\partial \theta_i}
     + \dfrac{\partial^2 L_n(T, \overline{\theta}_{n,i})}{ \partial \theta \partial \theta_i}(  \theta  - \theta_0)
 \end{equation}
   where $ [a,b]= \{ \lambda a + (1-\lambda)b ~ ; ~ \lambda \in [0,1]   \} $.
 Denote  $ \overline{F}_n(T, \theta)= -2\big(
   \dfrac{1}{card(T)} \dfrac{\partial^2 L_n(T, \overline{\theta}_{n,i})}{ \partial \theta \partial \theta_i}\big)_{1\leq i \leq
   d}$. Then, (\ref{Taylor_simple}) implies,
 \begin{equation}\label{Taylor_mat}
   Card(T)\overline{F}_n(T,\theta) ( \theta -\theta_0) = -2 \big( \dfrac{\partial L_n(T, \theta)}{\partial \theta} - \dfrac{\partial L_n(T,\theta_0)}{\partial \theta}  \big).
  \end{equation}
 Similarly, for any $\theta \in \Theta  $ we can find a matrix $ \widetilde{F}_n(T,\theta)$ such that

\begin{equation}\label{Taylor_mat_L_appro_1}
   Card(T)\widetilde{F}_n(T,\theta) ( \theta -\theta_0) = -2 \big( \dfrac{\partial \widehat{L}_n(T, \theta)}{\partial \theta} - \dfrac{\partial \widehat{L}_n(T,\theta_0)}{\partial \theta}  \big).
  \end{equation}
 With $ \theta = \widehat{\theta}_n(T) $ in (\ref{Taylor_mat_L_appro_1})
 and using the fact that $ \partial \widehat{L}_n(T,\widehat{\theta}_n(T))/ \partial \theta = 0 $  (because $ \widehat{\theta}_n(T) $ is a local extremum of $
 \widehat{L}_n(T,\cdot)$), it comes
 \begin{equation}\label{Taylor_mat_L_appro_2}
   Card(T) \widetilde{F}_n(T, \widehat{\theta}_n(T))( \widehat{\theta}_n(T)-\theta_0) = 2 \dfrac{\partial \widehat{L}_n(T,\theta_0)}{\partial \theta}.
 \end{equation}
 \begin{rem} \label{rem_conv_vers_F}
  If $Card(T) \limiten \infty  $ and  $ \theta = \theta(n)  \limiten \theta_0   $, then $  \overline{F}_n(T, \theta) \limitepsn F $ and
   $  \widetilde{F}_n(T, \theta)  \limitepsn   F $
  (see \cite{Bardet2009} and \cite{Bardet2010}). In particular,
  if $Card(T) \limiten  \infty $     , then $  \overline{F}_n(T, \widehat{\theta}_n(T))   \limitepsn   F $ and
    $ \widetilde{F}_n(T, \widehat{\theta}_n(T))  \limitepsn   F.$
 \end{rem}

\begin{lem} \label{lem_cov_F_til_theta_hat}
 Under assumptions of Theorem \ref{theo1}
 $$   \dfrac{1}{\sqrt{n}} ~ \underset{ k \in \Pi_n}{\mbox{max}} ~ \big\| k \big(\widetilde{F}_n(T_k,\widehat{\theta}_n(T_k))-F \big)(\widehat{\theta}_n(T_k)-\theta_0) \big\|=o_P(1).   $$

 \end{lem}
\begin{dem}
  For $k \in \Pi_n$, we know that $ \sqrt{k}( \widehat{\theta}_n(T_k)) - \theta_0 ) $
  converges in distribution to the Gaussian law as $ n \longrightarrow \infty $ (see Theorem 2 of \cite{Bardet2009}).
  Therefore,
  $ \underset{ k \in \Pi_n}{\mbox{max}} ~ \big\|\sqrt{k}(\widehat{\theta}_n(T_k)-\theta_0) \big\|=O_P(1) .$
  Remark \ref{rem_conv_vers_F}  implies that
  $ \underset{ k \in \Pi_n}{\mbox{max}} ~ \big\| \widetilde{F}_n(T_k,\widehat{\theta}_n(T_k))  -F \big\|=o(1) $ a.s.
  Thus
    \begin{align}
  \nonumber \dfrac{1}{\sqrt{n}} ~ \underset{ k \in \Pi_n}{\mbox{max}} ~ \big\| k \big(\widetilde{F}_n(T_k,\widehat{\theta}_n(T_k))-F \big)(\widehat{\theta}_n(T_k)-\theta_0) \big\|
   &\leq  \underset{ k \in \Pi_n}{\mbox{max}} ~ \big\| \widetilde{F}_n(T_k,\widehat{\theta}_n(T_k))  -F   \big\|
   \times \underset{ k \in \Pi_n}{\mbox{max}} ~ \big\|\sqrt{k}(\widehat{\theta}_n(T_k)-\theta_0) \big\|  \\
  \nonumber  & = o(1)O_P(1) ~ \text{a.s.} \\
  \nonumber  & = o_P(1).
  \end{align}
 \end{dem}
  Under assumptions of Theorem \ref{theo1}, the matrix $G$ is
  invertible. Denote $ \Sigma = FG^{-1}F $ \\
  $$  Q^{(1)}_n:= \underset{ k \in \Pi_n } {\mbox{max}}Q^{(1)}_{n,k} ~~ ~
 \text{where} ~~  Q^{(1)}_{n,k}:= \dfrac{k^2}{n} \big( \widehat{\theta}_n(T_k)-  \widehat{\theta}_n(T_n) \big)' \Sigma \big( \widehat{\theta}_n(T_k)-  \widehat{\theta}_n(T_n) \big) ~ ~ \text{and} ~ ~ ~ $$
   $$  Q^{(2)}_n:= \underset{ k \in \Pi_n } {\mbox{max}}Q^{(2)}_{n,k} ~~ ~
 \text{where} ~~  Q^{(2)}_{n,k}:= \dfrac{(n-k)^2}{n} \big( \widehat{\theta}_n(\overline{T}_k)-  \widehat{\theta}_n(T_n) \big)' \Sigma \big( \widehat{\theta}_n(\overline{T}_k)-  \widehat{\theta}_n(T_n) \big).$$

\begin{lem} \label{lem_Q_hat_Q}
 Under assumptions of Theorem \ref{theo1}
 $$  \underset{ k \in \Pi_n}{\mbox{max}} ~ \big| \widehat{Q}^{(j)}_{n,k} - Q^{(j)}_{n,k}  \big|=o_P(1) ~ ~ \text{for} ~ ~ j=1,2.   $$
 \end{lem}
 \begin{dem}
  The proof is provided for $j=1$, proceed similarly for $j=2$. For any $ k\in \Pi_n$, we have
  \begin{align}
  \nonumber \big| \widehat{Q}^{(1)}_{n,k} - Q^{(1)}_{n,k}  \big|  &\le \dfrac{k^2}{n} \|\widehat{\theta}_n(T_k)-  \widehat{\theta}_n(T_n)\|^2 \| \widehat{\Sigma}_{n,k}-\Sigma  \|  \\
  \nonumber  &\le  2 \dfrac{k^2}{n} \big( \|\widehat{\theta}_n(T_k)-  \theta_0\|^2 + \|\widehat{\theta}_n(T_n)-  \theta_0\|^2 \big) \| \widehat{\Sigma}_{n,k}-\Sigma  \| \\
  \label{lem_Q_hat_Q_eq1}    &\le 2 \big( \| \sqrt{k}(\widehat{\theta}_n(T_k)-  \theta_0)\|^2 + \|\sqrt{n}(\widehat{\theta}_n(T_n)-  \theta_0)\|^2 \big) \| \widehat{\Sigma}_{n,k}-\Sigma  \| .
\end{align}
 Since $k \in \Pi_n$, $k, n-k \longrightarrow \infty  $ as $ n \longrightarrow \infty $.
 Therefore, $ \sqrt{k}(\widehat{\theta}_n(T_k)-  \theta_0) = O_P(1)$ as $  n \longrightarrow \infty$,   $\sqrt{n}(\widehat{\theta}_n(T_n)-  \theta_0) = O_P(1)$,
   $ \widehat{F}_n(T_k) \limitepsn F $,  $ \widehat{F}_n(\overline{T}_k) \limitepsn F$,
  $ \widehat{G}_n(T_k) \limitepsn G $   and $ \widehat{G}_n(\overline{T}_k) \limitepsn G$ which is invertible.
 Thus, for n large enough, $ \widehat{G}_n(T_k)$ and  $  \widehat{G}_n(\overline{T}_k)$
 are invertible. It follows that as $ n \longrightarrow \infty$,
 \begin{align*}
    \| \widehat{\Sigma}_{n,k}-\Sigma  \| & = \big\|   \dfrac{k}{n}\widehat{F}_n(T_k)\widehat{G}_n(T_k)^{-1}\widehat{F}_n(T_k)
      +  \dfrac{n-k}{n}\widehat{F}_n(\overline{T}_k)\widehat{G}_n(\overline{T}_k)^{-1}\widehat{F}_n(\overline{T}_k)  -  FG^{-1}F \big\|\\
    & =  \big \|  \dfrac{k}{n}\big( \dfrac{k}{n}\widehat{F}_n(T_k)\widehat{G}_n(T_k)^{-1}\widehat{F}_n(T_k)- FG^{-1}F \big) +
        \dfrac{n-k}{n}\big( \widehat{F}_n(\overline{T}_k)\widehat{G}_n(\overline{T}_k)^{-1}\widehat{F}_n(\overline{T}_k)- FG^{-1}F \big)  \big\|\\
    & \leq  \| \widehat{F}_n(T_k)\widehat{G}_n(T_k)^{-1}\widehat{F}_n(T_k)- FG^{-1}F \|
         +  \|\widehat{F}_n(\overline{T}_k)\widehat{G}_n(\overline{T}_k)^{-1}\widehat{F}_n(\overline{T}_k)- FG^{-1}F  \|=  o(1) ~ \text{a.s.}
\end{align*}
 Therefore, (\ref{lem_Q_hat_Q_eq1}) implies
  $ \underset{ k \in \Pi_n}{\mbox{max}} ~ \big| \widehat{Q}^{(1)}_{n,k} - Q^{(1)}_{n,k}  \big| =o_P(1).$
 \end{dem}
\begin{lem} \label{lem_cov_derive}
 Under assumptions of Theorem \ref{theo1}
 $$ \dfrac{-2}{\sqrt{n}} \dfrac{\partial L_n (T_{[n\tau]},\theta_0)}{\partial \theta} \overset{ \mathcal{D}}{\underset{} \longrightarrow} ~~  W_G(\tau)
    ~~ \text{ in} ~~ D([0,1], \R^d )  $$
 where $W_G$ is a d-dimensional Gaussian process with zero mean and
 covariance matrix $\text{min}(\tau,s)G$.
 \end{lem}
\begin{dem}
 Recall that
 $ -2\dfrac{\partial L_n (T_{[n\tau]},\theta_0)}{\partial \theta}
  = \sum\limits_{t=1}^{[n\tau]} \dfrac{\partial q_{t} (\theta_0)}{\partial \theta}  .$
 Denote $ \mathcal{F}_t = \sigma (X_{t-1},\cdots) $. Since $X$ is
 stationary and ergodic, it is the same for the process  $(\dfrac{\partial q_{t} (\theta_0)}{\partial \theta})_{t\in \Z}$.
 Moreover,  $(\dfrac{\partial q_{t} (\theta_0)}{\partial \theta}, \mathcal{F}_t)$
 is a square integrable martingale difference  process (see \cite{Bardet2009}) with covariance matrix G.
 Then, the result follow by using Theorem 23.1 Billingsley (1968) (see \cite{Billingsley1968} page 206).
\end{dem}
\begin{lem} \label{lem_cov_derive_sum}
 Under assumptions of Theorem \ref{theo1}
 $$ \dfrac{-2}{\sqrt{n}} G^{-1/2} \Big( \dfrac{\partial L_n(T_{[n\tau]},\theta_0)}{\partial \theta}  -
    \dfrac{[n\tau]}{n}\dfrac{\partial L_n(T_n,\theta_0)}{\partial \theta} \Big )  \overset{ \mathcal{D}}{\underset{} \longrightarrow} ~~  W_d(\tau)
    ~~ \text{ in} ~~ D([0,1], \R^d )  $$
 where $W_d$ is a d-dimensional Brownian bridge.
 \end{lem}
\vspace{-0.5cm}
\begin{dem}
By Lemma \ref{lem_cov_derive}, it comes
$$ \dfrac{-2}{\sqrt{n}}\Big( \dfrac{\partial L_n(T_{[n\tau]},\theta_0)}{\partial \theta}  -
    \dfrac{[n\tau]}{n}\dfrac{\partial L_n(T_n,\theta_0)}{\partial \theta} \Big ) \overset{ \mathcal{D}}{\underset{} \longrightarrow} ~~
    W_G(\tau)-\tau W_G(1) ~~ \text{ in} ~~ D([0,1], \R^d ) .  $$
 Since the covariance matrix of the process $ \{  W_G(\tau)-\tau W_G(1) , ~ 0 \leq \tau \leq 1  \} $ is $ ( \text{min}(\tau , s) - \tau s )G $,
 the covariance matrix of the process $ \{  G^{-1/2}(W_G(\tau)-\tau W_G(1)) , ~ 0 \leq \tau \leq 1  \}$
  is $ ( \text{min}(\tau , s) - \tau s )I_d $ (where $I_d$ is the d-dimensional  identity matrix).
  Therefore, the process is equal (in distribution) to a
  d-dimensional Brownian bridge and the result follows.
\end{dem}
\begin{lem} \label{lem_cov_loi_Rd}
 Under assumptions of Theorem \ref{theo1}
$$ \dfrac{-2}{\sqrt{n}} G^{-1/2} \dfrac{\partial \widehat{L}_n ( T_{[n\tau]} , \widehat{\theta}_n(T_n) )}{\partial \theta}
      \overset{ \mathcal{D}}{\underset{} \longrightarrow} ~~  W_d(\tau)
    ~~ \text{ in} ~~ D([0,1], \R^d ) . $$

 \end{lem}
\begin{dem}
 From \cite{Bardet2009}, we have
  $ \dfrac{1}{\sqrt{n}} \big\| \dfrac{\partial L_n(T_n,\cdot)}{\partial \theta} - \dfrac{\partial \widehat{L}_n(T_n,\cdot)}{\partial \theta} \big\|_{\Theta}=o_P(1)$.
  This implies,
  \begin{equation}\label{conv_L_hat_L}
    \dfrac{1}{\sqrt{n}} ~ \underset{ k \in \Pi_n}{\mbox{max}} ~ \big\| \dfrac{\partial L_n(T_k,\cdot)}{\partial \theta}
   - \dfrac{\partial \widehat{L}_n(T_k,\cdot)}{\partial \theta}\big\|_{\Theta}=o_P(1).
   \end{equation}
 Let $k\in \Pi_n.$  Applying (\ref{Taylor_mat}) with $ T = T_k $ and $  \theta = \widehat{\theta}_n(T_n) $, we have
  $$   k\overline{F}_n(T_k,\widehat{\theta}_n(T_n)) ( \widehat{\theta}_n(T_n) -\theta_0) = -2 \big( \dfrac{\partial L_n(T_k, \widehat{\theta}_n(T_n))}{\partial \theta} -
    \dfrac{\partial L_n(T_k,\theta_0)}{\partial \theta}  \big)  .$$
 By plugging it in (\ref{conv_L_hat_L}), we have
\begin{equation}\label{lem_cov_loi_Rd_eq2}
    \dfrac{1}{\sqrt{n}} ~ \underset{ k \in \Pi_n}{\mbox{max}} ~ \big\|  \dfrac{\partial \widehat{L}_n(T_k,\widehat{\theta}_n(T_n))}{\partial \theta}
   -  \dfrac{\partial L_n(T_k,\theta_0)}{\partial \theta} + \dfrac{1}{2}k\overline{F}_n(T_k,\widehat{\theta}_n(T_n))(\widehat{\theta}_n(T_n) - \theta_0) \big\|=o_P(1).
   \end{equation}
 But, by Remark \ref{rem_conv_vers_F}, it comes that
\begin{multline*}
 \dfrac{1}{\sqrt{n}} ~ \underset{ k \in \Pi_n}{\mbox{max}} ~
  \big\|  k \big(  \overline{F}_n(T_k,\widehat{\theta}_n(T_n)) - \overline{F}_n(T_n,\widehat{\theta}_n(T_n)) \big)(\widehat{\theta}_n(T_n) - \theta_0) \big\|\\
  \hspace{4.65cm}   \leq \dfrac{1}{\sqrt{n}} ~ \underset{ k \in \Pi_n}{\mbox{max}} ~ \big\|  k \big(  \overline{F}_n(T_k,\widehat{\theta}_n(T_n)) -
   \overline{F}_n(T_n,\widehat{\theta}_n(T_n)) \big) \big\| \times \| \sqrt{n}(\widehat{\theta}_n(T_n) - \theta_0)\|\\
     = o(1)O_P(1) ~ ~ \text{a.s.} \hspace{3.9cm} \\
    = o_P(1). \hspace{9.68cm}
\end{multline*}
 Thus, (\ref{lem_cov_loi_Rd_eq2}) becomes
\begin{equation}\label{lem_cov_loi_Rd_eq3}
    \dfrac{1}{\sqrt{n}} ~ \underset{ k \in \Pi_n}{\mbox{max}} ~ \big\|  \dfrac{\partial \widehat{L}_n(T_k,\widehat{\theta}_n(T_n))}{\partial \theta}
   -  \dfrac{\partial L_n(T_k,\theta_0)}{\partial \theta} + \dfrac{1}{2}k\overline{F}_n(T_n,\widehat{\theta}_n(T_n))(\widehat{\theta}_n(T_n) - \theta_0) \big\|=o_P(1).
   \end{equation}
 Applying (\ref{Taylor_mat}) with $ T = T_n $ , $  \theta = \widehat{\theta}_n(T_n) $,
 and using $  (1/\sqrt{n})(\partial L_n(T_n,\widehat{\theta}_n(T_n)) /\partial \theta) = o_P(1) $     (see \cite{Bardet2009}), it follows
  \begin{equation}\label{lem_cov_loi_Rd_eq4b}
   \overline{F}_n(T_n,\widehat{\theta}_n(T_n)) ( \widehat{\theta}_n(T_n) -\theta_0) = \dfrac{2}{n}\dfrac{\partial L_n(T_n,\theta_0)}{\partial \theta}
                       +  o_P(\dfrac{1}{\sqrt{n}}) .
    \end{equation}
 Therefore, (\ref{lem_cov_loi_Rd_eq3}) becomes
 \begin{equation}\label{lem_cov_loi_Rd_eq4}
    \dfrac{1}{\sqrt{n}} ~ \underset{ k \in \Pi_n}{\mbox{max}} ~ \big\|  \dfrac{\partial \widehat{L}_n(T_k,\widehat{\theta}_n(T_n))}{\partial \theta}
   -  \dfrac{\partial L_n(T_k,\theta_0)}{\partial \theta} +  \dfrac{k}{n}\dfrac{\partial L_n(T_n,\theta_0)}{\partial \theta}\big\|=o_P(1).
   \end{equation}
 Now, let $0<\tau <1$, for large value of $n$, we have $  [\tau n] \in \Pi_n$;
 write
 \begin{multline*}
  \dfrac{-2}{\sqrt{n}} G^{-1/2} \dfrac{\partial \widehat{L}_n (T_{[n\tau]},\widehat{\theta}_n(T_n))}{\partial \theta}
     = \dfrac{-2}{\sqrt{n}} G^{-1/2} \Big[ \dfrac{\partial \widehat{L}_n (T_{[n\tau]},\widehat{\theta}_n(T_n))}{\partial \theta}
     -   \big( \dfrac{\partial L_n(T_{[n\tau]}, \theta_0)}{\partial \theta}  - \dfrac{[n\tau]}{n}\dfrac{\partial L_n (T_n,\theta_0)}{\partial \theta} \big )\\
    +  \big( \dfrac{\partial L_n(T_{[n\tau]}, \theta_0)}{\partial \theta}  -\dfrac{[n\tau]}{n}\dfrac{\partial L_n (T_n,\theta_0)}{\partial
     \theta} \big ) \Big ]
 \end{multline*}
 and the result follows by using  (\ref{lem_cov_loi_Rd_eq4}) and  Lemma \ref{lem_cov_derive_sum}.
\end{dem}

  \subsection{Proof of Theorem \ref{theo1}   and Theorem \ref{theo2}}
  ~

  \noindent {\bf Proof of Theorem \ref{theo1} .} \\
   We give the proof for $j=1$, proceed similarly for $j=2$. By
   Lemma \ref{lem_Q_hat_Q}, Theorem \ref{theo1} is established if
   $ Q^{(1)}_n   \limiteloin  \underset{ 0 \leq \tau \leq 1 } {\mbox{sup}} \| W_d(\tau) \|^2.$
   Using (\ref{conv_L_hat_L}), (\ref{Taylor_mat_L_appro_2}) with
   $T=T_k$ and Lemma \ref{lem_cov_F_til_theta_hat} it follows
    \begin{multline}
  \dfrac{1}{\sqrt{n}} ~ \underset{ k \in \Pi_n}{\mbox{max}} ~ \big\| \dfrac{\partial L_n(T_k,\theta_0)}{\partial \theta}
   - \dfrac{1}{2}k F~(\widehat{\theta}_n(T_k)-  \theta_0)  \big\|
   =  \dfrac{1}{\sqrt{n}} ~ \underset{ k \in \Pi_n}{\mbox{max}} ~ \big\| \dfrac{\partial \widehat{L}_n(T_k,\theta_0)}{\partial \theta}
   - \dfrac{1}{2}k F~(\widehat{\theta}_n(T_k)-  \theta_0)  \big\| + o_P(1) \\
   \hspace{1.67cm} = \dfrac{1}{\sqrt{n}} ~ \underset{ k \in \Pi_n}{\mbox{max}} ~ \big\| \dfrac{1}{2}k\widetilde{F}_n(T_k,\widehat{\theta}_n(T_k))(\widehat{\theta}_n(T_n)-  \theta_0)
   - \dfrac{1}{2}k F~(\widehat{\theta}_n(T_k)-  \theta_0)  \big\| +o_P(1)\\
   \label{theo1_eq1}  =  \dfrac{1}{\sqrt{n}} ~ \underset{ k \in \Pi_n}{\mbox{max}} ~ \big\| \dfrac{1}{2}k \big(\widetilde{F}_n(T_k,\widehat{\theta}_n(T_k))-F \big)
                      (\widehat{\theta}_n(T_n)-  \theta_0) \big\| +o_P(1) = o_P(1). \hspace{1.0cm}
 \end{multline}
 Using (\ref{lem_cov_loi_Rd_eq4}) and  \ref{theo1_eq1}, we have
   \begin{multline}
  \dfrac{1}{\sqrt{n}} ~ \underset{ k \in \Pi_n}{\mbox{max}} ~ \big\| \dfrac{\partial L_n(T_k, \widehat{\theta}_n(T_n))}{\partial \theta}
   - \dfrac{1}{2}k F~(\widehat{\theta}_n(T_k)-  \widehat{\theta}_n(T_n)) \big\|\\
   \hspace{3.2cm} = \dfrac{1}{\sqrt{n}} ~ \underset{ k \in \Pi_n}{\mbox{max}} ~ \big\|  \dfrac{\partial L_n(T_k,\theta_0)}{\partial \theta} -
     \dfrac{k}{n}\dfrac{\partial L_n(T_n,\theta_0)}{\partial \theta} - \dfrac{1}{2}k F~(\widehat{\theta}_n(T_k)-  \widehat{\theta}_n(T_n)) \big\|+o_P(1) \\
   \hspace{3.6cm} =  \dfrac{1}{\sqrt{n}} ~ \underset{ k \in \Pi_n}{\mbox{max}} ~ \big\|  \dfrac{1}{2}k F~(\widehat{\theta}_n(T_k)-  \theta_0) -
     \dfrac{k}{n}\dfrac{\partial L_n(T_n,\theta_0)}{\partial \theta} - \dfrac{1}{2}k F~(\widehat{\theta}_n(T_k)-  \widehat{\theta}_n(T_n)) \big\|+o_P(1) \\
  \hspace{0.18cm}=  \dfrac{1}{\sqrt{n}} ~ \underset{ k \in \Pi_n}{\mbox{max}} ~ \big\|  \dfrac{1}{2}k F~(\widehat{\theta}_n(T_n)-  \theta_0) -
     \dfrac{k}{n}\dfrac{\partial L_n(T_n,\theta_0)}{\partial \theta} \big\|+o_P(1) \\
 \label{theo1_eq2}     \leq  \sqrt{n} \big\|  \dfrac{1}{2} F~(\widehat{\theta}_n(T_n)-  \theta_0)
         - \dfrac{1}{n}\dfrac{\partial L_n(T_n,\theta_0)}{\partial \theta}  \big\| + o_P(1). \hspace{3.34cm}
 \end{multline}
 Note that
  \begin{align*}
  \big\| \sqrt{n} (F- \overline{F}_n(T_n,\widehat{\theta}_n(T_n))) ~(\widehat{\theta}_n(T_n)-  \theta_0)  \big\|
   & \leq   \big\| F- \overline{F}_n(T_n,\widehat{\theta}_n(T_n)) \big\| ~\big\| \sqrt{n}(\widehat{\theta}_n(T_n)-  \theta_0)   \big\| \\
    & = o(1) O_P(1) ~ \text{ a.s.} \\
    & = o_P(1).
\end{align*}
  By plugging it in (\ref{theo1_eq2})  and applying (\ref{Taylor_mat}) with $T=T_n$ and $ \theta = \widehat{\theta}_n(T_n)$, we have
\begin{align}
 \nonumber   \dfrac{1}{\sqrt{n}}\underset{ k \in \Pi_n}{\mbox{max}}\big\| \dfrac{\partial L_n(T_k, \widehat{\theta}_n(T_n))}{\partial \theta}
   - \dfrac{1}{2}k F~(\widehat{\theta}_n(T_k)-  \widehat{\theta}_n(T_n)) \big\|
   & \leq  \sqrt{n} \big\|  \dfrac{1}{2} \overline{F}_n(T_n,\widehat{\theta}_n(T_n))(\widehat{\theta}_n(T_n)-  \theta_0)
         - \dfrac{1}{n}\dfrac{\partial L_n(T_n,\theta_0)}{\partial \theta}  \big\| \\
   \label{theo1_eq3}  &   \hspace{5.0cm}    +   o_P(1).
\end{align}
 Therefore, using (\ref{lem_cov_loi_Rd_eq4b}), (\ref{theo1_eq3}) implies
 \begin{equation}\label{theo1_eq4}
  \dfrac{1}{\sqrt{n}}\underset{ k \in \Pi_n}{\mbox{max}}\big\| \dfrac{\partial L_n(T_k, \widehat{\theta}_n(T_n))}{\partial \theta}
   - \dfrac{1}{2}k F~(\widehat{\theta}_n(T_k)-  \widehat{\theta}_n(T_n))  \big\| = o_P(1).
      \end{equation}
  Now, let $0<\tau <1$, for large value of $n$, we have $  [\tau n] \in \Pi_n$; write
\begin{multline*}
  \dfrac{-2}{\sqrt{n}}G^{-1/2} \dfrac{\partial \widehat{L}_n(T_{[n \tau]},\widehat{\theta}_n(T_n))}{\partial \theta}
    =  - \dfrac{[n \tau]}{\sqrt{n}}G^{-1/2} F ( \widehat{\theta}_n(T_{[n \tau]}) - \widehat{\theta}_n(T_n))\\
   - 2G^{-1/2} \dfrac{1}{\sqrt{n}} \Big[  \dfrac{\partial \widehat{L}_n(T_{[n \tau]},\widehat{\theta}_n(T_n))}{\partial \theta}
    -  \dfrac{1}{2} [n \tau] F ( \widehat{\theta}_n(T_{[n \tau]}) - \widehat{\theta}_n(T_n) ) \Big].
 \end{multline*}
 Therefore, using (\ref{theo1_eq4}) we have
 $$  - \dfrac{[n \tau]}{\sqrt{n}}G^{-1/2} F ( \widehat{\theta}_n(T_{[n \tau]}) - \widehat{\theta}_n(T_n))
    =   \dfrac{-2}{\sqrt{n}}G^{-1/2} \dfrac{\partial \widehat{L}_n(T_{[n \tau]},\widehat{\theta}_n(T_n))}{\partial \theta} + o_P(1) $$
 and the result follows by using  Lemma \ref{lem_cov_loi_Rd}. \Box

~\\

 \noindent {\bf Proof of Theorem \ref{theo2} .} \\
 Let $ \tau^* \in (0,1) $ the true value of break. Denote $ k^*=[n\tau^*]$.
 For $n$ large enough , $ k^* \in \Pi_n$. Therefore, we have for $ j=1,2 $,
  $ ~ \widehat{Q}^{(j)}_n = \underset{ k \in \Pi_n}{\mbox{max}} \widehat{Q}^{(j)}_{n,k} ~\geq \widehat{Q}^{(j)}_{n,k^*} .$
 Thus, it follows that
 \begin{equation}\label{theo2_eq1}
  \widehat{Q}_n = \text{max} ( \widehat{Q}^{(1)}_n , \widehat{Q}^{(2)}_n ) \geq
  \text{max} ( \widehat{Q}^{(1)}_{n,k^*} ,\widehat{Q}^{(2)}_{n,k^*}).
 \end{equation}
 Since $ \theta^*_1, \theta^*_2 \in \overset{\circ}{\Theta}(4)$, it comes from \cite{Bardet2009} that the model
  $ \mathcal{M}_{\Z}(M_{\theta^*_1},f_{\theta^*_1})  $ and $ \mathcal{M}_{\Z}(M_{\theta^*_2},f_{\theta^*_2})  $
 have a 4-order stationary solution which we denote $( X_{t,j})_{t \in \Z} $ ~for $j=1,2$.\\
 For $j=1,2$~ denote for any ~$t \in \Z$,
 $
 q_{t,j}(\theta):=(X_{t, j}-f_\theta^{t,j})^2/(h^{t,j}_{\theta})+\log (h^{t,j}_\theta) $
 with  $f_\theta^{t,j}:=f_\theta(X_{t-1,j},X_{t-2,j},\ldots)$, $h^{t,j}_{\theta}:=(M_\theta^{t,j})^2$
  where $M_\theta^{t,j}:=M_\theta(X_{t-1,j},X_{t-2,j},\ldots)$.
  Also denote for $j = 1, 2$
  $$ F^{(j)} = \E[   \dfrac{\partial^2 q_{0,j}(\theta^*_j)}{\partial \theta \partial \theta '}  ]   ~~ \text{and} ~~
    G^{(j)} = \E\big[   \big(   \dfrac{\partial q_{0,j}(\theta^*_j)}{\partial \theta} \big) \big(   \dfrac{\partial q_{0,j}(\theta^*_j)}{\partial \theta} \big)'   \big]  .$$
 For $j=1,2$, Lemma \ref{lem_inv_G} implies that the matrix  $G^{(j)}$ is symmetric positive definite and Corollary 5.1 of
 \cite{Bardet2010} implies $  \widehat{G}_n(T_{k^*}) \limitepsn G^{(1)} $ and $  \widehat{G}_n(\overline{T}_{k^*})\limitepsn G^{(2)}$.
   Lemma 4 of \cite{Bardet2009} implies  $  \widehat{F}_n(T_{k^*}) \limitepsn F^{(1)} $ and $  \widehat{F}_n(\overline{T}_{k^*})\limitepsn F^{(2)}$.
 Therefore, it follows that
 \begin{multline}
   \widehat{\Sigma}_{n,k^*} := \dfrac{k^*}{n}\widehat{F}_n(T_{k^*})\widehat{G}_n(T_{k^*})^{-1}\widehat{F}_n(T_{k^*})\text{\large{\1}}_{\det(\widehat{G}_n(T_{k^*}))\neq 0 }
 + \dfrac{n-k^*}{n}\widehat{F}_n(\overline{T}_{k^*})\widehat{G}_n(\overline{T}_{k^*})^{-1}\widehat{F}_n(\overline{T}_{k^*})\text{\large{\1}}_{\det(\widehat{G}_n(\overline{T}_{k^*}))\neq 0 }\\
  \label{theo2_eq2} \limitepsn  \tau^* F^{(1)}(G^{(1)})^{-1}F^{(1)} + (1-\tau^*)F^{(2)}(G^{(2)})^{-1}F^{(2)}.
 \end{multline}
 Denote $ \Sigma = \tau^* F^{(1)}(G^{(1)})^{-1}F^{(1)} + (1-\tau^*)F^{(2)}(G^{(2)})^{-1}F^{(2)}$.
 It is easy to see that $\Sigma$ is a symmetric positive definite matrix.\\
 For all $\rho>0$ and $ \theta \in \Theta  $, denote $ B_o(\theta,\rho) $
 (rep. $B_c(\theta,\rho)$ ) the open (resp. closed) ball centered at
 $\theta$ of radius $ \rho$ in $\Theta$. i.e.
 $$   B_o(\theta,\rho) = \{ x \in \Theta ~ ; ~  \| \theta - x \| < \rho   \}  ~~ \text{and} ~~
   B_c(\theta,\rho) = \{ x \in \Theta ~ ; ~  \| \theta - x \| \leq \rho   \}  .$$
 For $ A\subset \Theta $, we denote $ A^c = \{ x \in \Theta ~~ ; ~~ x \notin A   \}.$\\
 Since $ \theta^*_1 \neq  \theta^*_2 $ and $ \theta^*_1 ,  \theta^*_2 \in \overset{\circ}{\Theta}(4) \subset \overset{\circ}{\Theta} $,
 then there exist $ \rho_1 >0$ and  $ \rho_2 >0$ such as $ B_o(\theta^*_1,\rho_1) \cap B_o(\theta^*_2,\rho_2) = \emptyset $.\\
  For all $n \in \N$, denote
   $$  \delta^{(j)}_n = \underset{x \in B_c(\theta^*_j,\rho_j/2) ;~ y \in B^c_o( \theta^*_j,\rho_j )  }{\mbox{inf}}
     \big(  (x-y)' \widehat{\Sigma}_{n,k^*} (x-y)    \big)  ~~ \text{for}~~ j=1,2.$$
  Also denote $  \delta^{(j)} = \underset{x \in B_c(\theta^*_j,\rho_j/2) ;~ y \in B^c_o( \theta^*_j,\rho_j )  }{\mbox{inf}}
     \big(  (x-y)'\Sigma (x-y)    \big) $.  It is easy to
     see that  $  \delta^{(j)} > 0$   for $ j=1,2$. \\ 
 Using (\ref{theo2_eq2}), we have
 \begin{equation}\label{theo2_eq3}
   \delta^{(j)}_n \limitepsn  \delta^{(j)} ~~ \text{for} ~~ j=1,2.
 \end{equation}
   From \cite{Bardet2009}  and \cite{Bardet2010}, we have $ \widehat{\theta}_n(T_{k^*}) \limitepsn  \theta^*_1 $
   and $ \widehat{\theta}_n(\overline{T}_{k^*}) \limitepsn  \theta^*_2$.
  Therefore, for $n$ large enough, $ \widehat{\theta}_n(T_{k^*}) \in B_o(\theta^*_1, \rho_1/2  ) $
   ~ and ~ $ \widehat{\theta}_n(\overline{T}_{k^*}) \in B_o(\theta^*_2, \rho_2/2  )$.
  Thus, two situations may occur
  \begin{itemize}
    \item if  $\widehat{\theta}_n(T_n) \in  B_o(\theta^*_2, \rho_2  ) $ i.e. $\widehat{\theta}_n(T_n) \in  B^c_o(\theta^*_1, \rho_1 )$
     then
  $ (\widehat{\theta}_n(T_{k^*}) - \widehat{\theta}_n(T_{n}))' \widehat{\Sigma}_{n,k^*}(\widehat{\theta}_n(T_{k^*}) - \widehat{\theta}_n(T_{n}))
    \geq \delta^{(1)}_n$. Therefore,
    $$ \widehat{Q}^{(1)}_{n,k^*}:= \dfrac{(k^*)^2}{n}(\widehat{\theta}_n(T_{k^*}) - \widehat{\theta}_n(T_{n}))' \widehat{\Sigma}_{n,k^*}(\widehat{\theta}_n(T_{k^*}) - \widehat{\theta}_n(T_{n}))
     \geq \dfrac{(k^*)^2}{n}\delta^{(1)}_n \simeq n (\tau^*)^2 \delta^{(1)}_n . $$

    \item else  $\widehat{\theta}_n(T_n) \in  B^c_o(\theta^*_2, \rho_2 )$  and
  $ (\widehat{\theta}_n(\overline{T}_{k^*}) - \widehat{\theta}_n(T_{n}))' \widehat{\Sigma}_{n,k^*}(\widehat{\theta}_n(\overline{T}_{k^*}) - \widehat{\theta}_n(T_{n}))
    \geq \delta^{(2)}_n$. Therefore,
 $$ \widehat{Q}^{(2)}_{n,k^*} = \dfrac{(n-k^*)^2}{n}(\widehat{\theta}_n(\overline{T}_{k^*}) - \widehat{\theta}_n(T_{n}))' \widehat{\Sigma}_{n,k^*}(\widehat{\theta}_n(\overline{T}_{k^*}) - \widehat{\theta}_n(T_{n})
     \geq \dfrac{(n-k^*)^2}{n}\delta^{(2)}_n \simeq n (1-\tau^*)^2 \delta^{(2)}_n   .$$
  \end{itemize}
  In all cases, we have
  $ \widehat{Q}_n \geq \text{max}( \widehat{Q}^{(1)}_{n,k^*} , \widehat{Q}^{(2)}_{n,k^*}  )  \geq \text{min}\big( n (\tau^*)^2 \delta^{(1)}_n , n (1-\tau^*)^2 \delta^{(2)}_n \big) .$

  Thus the result follows by using (\ref{theo2_eq3}). \Box

 \section*{ Acknowledgements}
  The author thanks Jean-Marc Bardet and Olivier Wintenberger for many  discussions which
  helped to improve this work.

\end{document}